\documentclass[12pt,myheadings]{article}
\usepackage{natbib}
\usepackage{authblk}
\usepackage{fullpage}
\usepackage[para]{footmisc}
\usepackage{amsfonts,amsmath,amssymb,mathtools,amsthm}
\usepackage{listings}
\usepackage{graphicx}
\usepackage{booktabs}
\usepackage{url}

\usepackage{breakurl}
\usepackage[breaklinks]{hyperref}

\def\ignore#1{{}}

\ignore{
\makeatletter
\DeclareRobustCommand{\cev}[1]{%
  {\mathpalette\do@cev{#1}}%
}
\newcommand{\do@cev}[2]{%
  \vbox{\offinterlineskip
    \sbox\z@{$\m@th#1 x$}%
    \ialign{##\cr
      \hidewidth\reflectbox{$\m@th#1\vec{}\mkern4mu$}\hidewidth\cr
      \noalign{\kern-\ht\z@}
      $\m@th#1#2$\cr
    }%
  }%
}
\makeatother
}
\def\bbbone{{\mathchoice {\rm 1\mskip-4mu l} {\rm 1\mskip-4mu l}
{\rm 1\mskip-4.5mu l} {\rm 1\mskip-5mu l}}}

\renewcommand{\binom}[2]{\genfrac{(}{)}{0pt}{}{#1}{#2}}

\def\bbbn{\mathbb{N}}
\def\ignore#1{}

\def\bbbone{{\mathchoice {\rm 1\mskip-4mu l} {\rm 1\mskip-4mu l}
{\rm 1\mskip-4.5mu l} {\rm 1\mskip-5mu l}}}

\def\BZ{\mathbb{Z}}

\def\bbbn{\mathbb{N}}
\def\ignore#1{}

\def\BP{\mathbb{P}}
\def\BE{\mathbb{E}}
\def\BN{\mathbb{N}}

\def\eq{\begin{equation}}
\def\en{\end{equation}}

\newcommand{\C}{\mathcal{C}}

\newtheorem{proof*}{Proof}
\newtheorem{corollary}{Corollary} 
\newtheorem{theorem}{Theorem}

\newtheorem{lemma}{Lemma}

 \makeatletter
   \def\@fnsymbol#1{\ensuremath{\ifcase#1\or 1\or
   2\or 3\or \|\or **\or 4
   \or 5 \else\@ctrerr\fi}}
    \makeatother
\def\BE{\mathbb{E}}
\def\BZ{\mathbb{Z}}
\def\BP{\mathbb{P}}
\def\bbbn{\mathbb{N}}
\def\BN{\mathbb{N}}
\def\ignore#1{{}}
\def\C{\mathcal{C}}
\def\bbbone{{\mathchoice {\rm 1\mskip-4mu l} {\rm 1\mskip-4mu l}
{\rm 1\mskip-4.5mu l} {\rm 1\mskip-5mu l}}}
\usepackage{graphicx}

\providecommand{\keywords}[1]
{
  \small	
 \noindent\textbf{Keywords:} #1
}

\providecommand{\msc}[1]
{
  \small	
 \noindent\textbf{MSC:} #1
}

\begin{document}

\pagestyle{plain}
\title{Another view of sequential sampling in the birth process with immigration\let\thefootnote\relax\footnotetext{\textit{E-mail addresses:} phd2120@columbia.edu,  aj2963@columbia.edu, st3193@columbia.edu}}

\author[a,c]{Poly H. da Silva
}
\affil[a]{{\footnotesize Columbia University, Department of Statistics, 1255 Amsterdam Avenue, New York, NY 10027, USA}}
\author[b,c]{Arash Jamshidpey
}
\affil[b]{{\footnotesize Columbia University, Department of Mathematics, 2990 Broadway, New York, NY 10027, USA}}
\author[a,c]{Simon Tavar\'e
}
\affil[c]{{\footnotesize Columbia University, Irving Institute for Cancer Dynamics, Schermerhorn Hall, Suite 601, 1190 Amsterdam Avenue, New York, NY 10027, USA}}
\date{\begin{footnotesize}
November 28, 2022
\end{footnotesize}}
\maketitle


\abstract{Models of counts-of-counts data have been extensively used in the biological sciences, for example in  cancer, population genetics, sampling theory and ecology. In this paper we explore properties of one model that is embedded into a continuous-time process and can describe the appearance of certain biological data such as covid DNA sequences in a database. More specifically, we consider an evolving model of counts-of-counts data that arises as the family size counts of samples taken sequentially from a Birth process with Immigration (BI).  Here, each family represents a type or species, and the family size counts represent the type or species frequency spectrum in the population. We study the correlation of $S(a,b)$ and $S(c,d)$, the number of families observed in two disjoint time intervals $(a,b)$ and $(c,d)$. We find the expected sample variance and its asymptotics for $p$ consecutive sequential samples $\mathbf{S}_p:=(S(t_0,t_1),\dots, S(t_{p-1},t_p))$, for any given $0=t_0<t_1<\dots<t_p$. By conditioning on the sizes of the samples, we provide a connection between $\mathbf{S}_p$ and $p$ sequential samples of sizes $n_1,n_2,\dots,n_p$,  drawn from a single run of a Chinese Restaurant Process. The properties of the latter were studied in~\cite{SJMT22}. We show how the continuous-time framework helps to make asymptotic calculations easier than its discrete-time counterpart. As an application, for a specific choice of $t_1,t_2,\dots, t_p$, we revisit Fisher's 1943 multi-sampling problem and give another explanation of what Fisher's model could have meant in the world of sequential samples drawn from a BI process.}



\msc{60C05, 60J25, 60J80, 92D10, 92D25, 92D40}

\maketitle

\section{Introduction}\label{sec1}

The work we describe here arose from an attempt to model the arrival of covid DNA sequences at the GISAID database \citep{GISAID21}. The sequences arrive sequentially through time, a first sequence, a second, and so on. One model might therefore be a discrete-time process in which sequences are labeled by their order of appearance, and the type of that sequence (for example, the sequence itself) is recorded for each of them. One standard model for such counts-of-counts data is the so-called Chinese Restaurant Process (CRP), which records the accumulation of different sequence types (referred to as families in what follows) as they arrive at the database: a newly arriving sequence is either a copy of an existing sequence, or a novel sequence, with a probability determined by a parameter $\theta$ and the arrival number of the sequence. \cite{SJMT22} describe such a model in some detail, and analyses, \textit{inter alia}, the behavior of the number of distinct families, $S_i$, seen in sequential samples of size $n_i$, $i = 1, 2, \ldots, p.$ 

In fact the samples arrive continuously through time, so that the samples are (for example) those obtained in week 1, week 2, and so on. This requires a continuous-time process as its generative model, and here we chose the Yule process with immigration. We note that this is \emph{not} intended to provide a mechanistic model for the way the covid sequences have evolved, but rather  a description of their arrival at the database. Thus we do not need to model deaths in our process.

This simplification leads to a model for which explicit results concerning the behavior of the number of families arising in $p$ sequential time intervals may be derived. The analysis relies on simple Marked Poisson Process arguments, for which various marking probabilities may be calculated explicitly.

The work in \cite{SJMT22} (see also \cite{dgk48a,dgk49}, \cite{eclw2007}, \cite{bt10}) was motivated by an interpretation of a much older problem of \cite{raf43}, who was concerned with finding the expected sample variance of (in our phraseology) $S_1$ and $S_2$, namely $\BE V_2 := \BE (S_1 - S_2)^2/2$, when the two sample sizes are equal, to $n$ say. Fisher surmised that when $n$ is large, $\BE V_2 \sim \theta \log 2$. The present results provide an elementary analysis of this problem in a more general setting, and provides another view of Fisher's calculations. 

\section{The Yule process with immigration}\label{bdi}
We begin with some well-known results for the Yule process with immigration (cf. \cite{st87}). Immigration events occur at the points of a Poisson process of rate $\theta$, each new immigrant initiating a family that grows according to a pure birth (Yule) process, $B(\cdot)$, with time scaled so that the birth rate is $\lambda = 1$. (This results in no loss of generality, since  for arbitrary $\lambda$ one replaces  $t$ by $\lambda t$, and immigration rate $\theta = \theta/\lambda$.)
The distribution of the number of members of a typical family that has grown for time $t$ from a single individual is geometric:
\begin{equation}\label{bprocess}
p_j(t) := \BP(B(t) = j \vert B(0) = 1) = e^{-t} (1-e^{-t})^{j-1}, j=1,2,\ldots
\end{equation}
We denote the transition probability $\BP(B(t) = k \vert B(0) = j)$ by $p_{jk}(t).$

The population size at time $t$ is denoted by $Z(t)$, and
the number of families of size $i$ at time $t$ is denoted by $C_i(t), i =1, 2, \ldots.$
It is well known (cf. \cite{km67, st87}) that for each $t$, the $C_i(t)$ are independent Poisson random variables with means given by
\begin{equation}\label{Ctmeans}
\BE C_i(t) = \theta (1-e^{-t})^i / i, i=1,2,\ldots.
\end{equation}

What is less clear is the joint behavior of the number of families observed in different time intervals. This note addresses aspects of this problem.

We make use of the probability generating function (pgf) of $B(t)$ which is given by
\begin{equation}\label{pgf}
\phi(t;s) := \BE s^{B(t)} = \sum_{j \geq 1} p_j(t) s^j = \frac{e^{-t} s}{1 - (1 - e^{-t})s},
\end{equation}
for $0 \leq s \leq 1$, 
whence it follows that for $0 \leq a \leq b \leq c$, we have
\begin{equation}\label{intresult}
\Phi(a,b,c;s) := \int_a^b \,\phi(c-u;s) du = \log\left( \frac{1 - (1 - e^{-(c-b)}) s}{  1 - (1 - e^{-(c-a)}) s}\right).
\end{equation}

\section{Families observable in a given period of time}\label{sect3}
If we were to watch the growth of a given family through time, there will be some intervals in which no new births are observed. In such an interval, all the members of the family were old: that is, were born before the start of the interval; otherwise, the family is composed of old members and new ones. It is the interplay between these two types that we uncover.

We say a family is \emph{observable} in $(a,b)$, for $0\leq a<b$, if it has at least one birth in $(a,b)$. Let $S(a,b)$ count the families observable in $(a,b)$. Of interest is the distribution of $S(a,b)$ and the joint distribution of $S(a,b)$ and $S(c,d)$ for $0\leq a<b\leq c<d$. In particular, we compute ${\rm Cov}(S(a,b),S(c,d))$. 

Letting $I_1=(a,b), I_2=(c,d)$, we notice that
\begin{eqnarray*}
S(I_1)&=&S(a,b)=K(I_1,I_2)+T(I_1\setminus I_2),\\
S(I_2)&=&S(c,d)=K(I_1,I_2)+T(I_2\setminus I_1),
\end{eqnarray*}
where $K(I_1,I_2)$ counts families observable in both $I_1$ and $I_2$, and $T(I_i\setminus I_j)$ counts families observable in $I_i$ but not in $I_j$, for $\{i\neq j\}=\{1,2\}$. In fact, $K(I_1,I_2)$ counts those points of the Poisson process $X$ on $(0,d)$, which are marked if observed in both $I_1$ and $I_2$. Similar marking applies to $T(I_1\setminus T_2)$ and $T(I_2\setminus T_1)$. From the Marking Theorem (cf. \cite{jfck93}),  $K(I_1,I_2)$, $T(I_1\setminus T_2)$ and $T(I_2\setminus T_1)$ are independent Poisson random variables. To obtain their expected values,  for $J=(a_0,b_0)$ and $I_i=(a_i,b_i)$ such that $b_{i-1}\leq a_i$ for $i=1,\dots,k$, we define random variables $V_J(I_1,\dots,I_k)=V_{a_0,b_0}(I_1,\dots,I_k)$ to count the number of families initiated in $J$ and not observable in $I_1,\dots,I_k$. We need the following lemma.

\begin{lemma}\label{lemma1}
For $0\leq a<b\leq c<d$, let $J=(a,b)$ and $I=(c,d)$. Then $V_J(I)$ is a Poisson r.v.  with mean
$$
\BE V_J(I) = \theta \log\left( \frac{e^d - e^c + e^b}{e^d - e^c + e^a}\right)
$$
\end{lemma}

\medskip
\noindent{\bf Proof. }
For $r\in\bbbn$, let $X_{J,r}(I)$ count all those families initiated in $J$, with exactly $r$ members at $c$, which are not observable in $I$. The probability that a family initiated at $x\in J$, is not observable in $I$, while it has exactly $r$ members at $c$ is given by
$$
p_{r}(c-x)e^{-r (d-c)}.
$$
Hence, from the Marking Theorem, for $r\in \bbbn$, $V_{J,r}(I)$ are independent Poisson r.v.s with expected value
$$
\BE V_{J,r}(I) = \theta \int_J p_{r}(c-x)e^{-r(d-c)}dx.
$$
As a consequence, $V_J(I)$ is a Poisson r.v. with parameter 
\begin{eqnarray*}
\BE V_J(I) &=& \sum_{r=1}^\infty \BE V_{J,r}(I) 
= \theta \int_a^b \phi(c-x;e^{-(d-c)}) dx \\
&=& \theta \Phi(a,b,c;e^{-(d-c)}) 
 =  \theta \log\left( \frac{e^d - e^c + e^b}{e^d - e^c + e^a}\right),
\end{eqnarray*}
using (\ref{intresult}) and simplifying.
\hfill\qed

We highlight one special case of this result for the case $J = (0,a), I = (a,b)$. Lemma~\ref{lemma1} gives 
$$
\BE V_J(I) = \theta b - \theta \log(e^b - e^a + 1)
$$
We can establish the following lemma in a similar way.
\begin{lemma}\label{lemma2}
For $0\leq a<b\leq c<d$, let $J=(0,a)$, $I_1 = (a,b)$ and $I_2=(c,d)$. Then $V_J(I_1,I_2)$ is a Poisson r.v.  with expected value
$$
\BE V_J(I_1,I_2) = \theta \log\left( \frac{e^d - e^c + e^b}{e^d - e^c + e^b - e^a +1}\right)
$$
\end{lemma}

\medskip
\noindent{\bf Proof. }
The probability that a family arrives in $J$ and is not observable in $I_1$ and $I_2$, while it has $r$ and $s$ members at $a$ and $c$, respectively, $1\leq r\leq s$, is given by
$$
p_r(a-x)e^{- r(b-a)}p_{rs}(c-b)e^{- s(d-c)}.
$$
Hence, from the Marking Theorem, for $b<c$, $V_J(I_1,I_2)$ is a Poisson r.v.  with
\begin{eqnarray*}
\BE V_J(I_1,I_2)&  = & \sum_{s\geq 1}\sum_{r\leq s}\theta\int_0^a p_{r}(a-x)e^{-r(b-a)} p_{rs}(c-b)e^{-s(d-c)}dx\\
& = & \theta\int_0^a \sum_{r\geq 1} p_{r}(a - x) e^{-r(b-a)} \sum_{s\geq r} p_{rs}(c-b) e^{-s(d-c)}\\
& = & \theta \int_0^a \sum_{r\geq 1} p_{r}(a - x) e^{-r(b-a)} [\phi(c-b;e^{-(d-c)})]^r \\
& = & \theta \Phi(0,a,a;e^{-(b-a)}\phi(c-b;e^{-(d-c)})),
\end{eqnarray*}
the second to last line coming from the fact that families evolve independently.
Using (\ref{intresult}) once more, and simplifying, we get
$$
\BE V_J(I_1,I_2)  = \theta \log\left( \frac{e^d - e^c + e^b}{e^d - e^c + e^b - e^a +1}\right),
$$
as required. The case $b = c$ reduces to $\BE V_J(I)$ for $J = (0,a), I = I_1 \cup I_2 = (a,d),$ given in Lemma \ref{lemma1}.\hfill\qed

To compute the expected values of $K(I,J), T(I\setminus J)$ and $T(J\setminus I)$, we also need to define the random variables $U_J(I_1,\dots, I_k)$ to count the families initiated in $J$ and observable in $I_1,\dots,I_k$.

\begin{lemma}\label{lemma3}
For $0\leq a<b\leq c<d$, let $J=(a,b)$ and $I=(c,d)$. Then $U_J(I)$ is a Poisson r.v. with parameter
$$
\BE U_J(I) = \theta \log\left( \frac{e^b(e^d - e^c + e^a)}{e^a(e^d - e^c + e^b)}\right)
$$
\end{lemma}

\medskip
\noindent{\bf Proof. }
The Marking Theorem shows that $U_J(I)$ and $V_J(I)$ are independent Poisson r.v.s. The result follows since $X(J) = U_J(I) + V_J(I)$.
\hfill\qed

\begin{corollary}\label{thm2}
For any $0\leq a<b$, $S(a,b)$ is a Poisson random variable with mean
$$
\BE S(a,b)= \theta \log(e^b - e^a +1).
$$
\end{corollary}

\medskip\noindent{\bf Proof. }
Letting $I=(a,b)$, the theorem follows from Lemma~\ref{lemma3} and the fact that 
$S(a,b)=U_{(0,a)}(I)+X(I)$. Note that for an interval $(0,b)$, $S(0,b)$ has a Poisson distribution with mean $\theta \log(e^b -e^0 +1) = \theta \log(e^b) = \theta b,$ as it must.
\hfill\qed

\begin{theorem}\label{thm1}
For $0\leq a<b\leq c<d$, let $I_1=(a,b)$ and $I_2=(c,d)$.  Then $K(I_1,I_2)$, $T(I_1\setminus I_2)$ and $T(I_2\setminus I_1)$ are independent Poisson random variables with means
\begin{equation*}
\begin{split}
&\BE K(I_1,I_2) = \theta \log\left( \frac{(e^b - e^a + 1)(e^d -e^c + 1)}{e^d - e^c + e^b -e^a + 1}\right), \\
&\BE T(I_1\setminus I_2)  = \theta\log\left(\frac{e^d - e^c + e^b -e^a + 1}{e^d -e^c + 1}
\right),\\
&\BE T(I_2\setminus I_1)=\theta\log\left(\frac{e^d - e^c + e^b -e^a + 1}{e^b -e^a + 1}
\right).\\
\end{split}
\end{equation*}
\end{theorem}

\medskip
\noindent{\bf Proof. }
Let $J=(0,a), J'=(b,c)$, $K(I_1,I_2)=U_{I_1}(I_2)+U_J(I_1,I_2)$, where 
$$
U_J(I_1,I_2)=X(J)-V_J(I_1)-V_J(I_2)+V_J(I_1,I_2).
$$
From Lemma~\ref{lemma1} and Lemma~\ref{lemma2}, we have
\begin{eqnarray*}
\BE U_J(I_1,I_2) & = & \theta (a-b) + \theta \log(e^b - e^a + 1) \\
& & - \theta \log(e^d - e^c + e^a) + \theta \log(e^d - e^c + 1)\\
& & + \theta \log(e^d - e^c + e^b) - \theta \log(e^d - e^c + e^b -e^a + 1)
\end{eqnarray*}
$\BE K(I_1,I_2)$ follows after using Lemma~\ref{lemma3} and simplifying.

On the other hand, 
\begin{eqnarray*}
T(I_1\setminus I_2)&=& V_{I_1}(I_2)+V_J(I_2)-V_J(I_1,I_2),\\[10pt]
T(I_2\setminus I_1)&=& V_J(I_1)-V_J(I_1,I_2)+U_{J'}(I_2)+X(I_2)\\
&=& V_J(I_1)-V_J(I_1,I_2)+X(J')-V_{J'}(I_2)+X(I_2),
\end{eqnarray*}
which, applying Lemma \ref{lemma1} and Lemma \ref{lemma2} and simplifying, gives the result. 
\hfill\qed

\begin{corollary}\label{thm3}
For $0\leq a< b\leq c<d$, we have
$$
{\rm Cov}(S(a,b),S(c,d))=
\theta \log\left(\frac{(e^b -e^a + 1) (e^d - e^c + 1)}{e^d -e^c + e^b -e^a + 1}\right).
$$
\end{corollary}

\medskip\noindent{\bf Proof. }
It is clear from Theorem~\ref{thm1} that
$${\rm Cov}(S(a,b),S(c,d))={\rm Var} K((a,b),(c,d)). $$
\hfill\qed

\section{Sampling in multiple intervals}

We consider sampling from $p$ intervals determined by the points  $t_0=0<t_1<\dots<t_p$. Let $\delta_i=t_i-t_{i-1}$ for $1\leq i\leq p$. For $i=1,\dots,p$, let $S_i := S(t_{i-1},t_i)$ be the number of families observable in the time interval $(t_{i-1},t_i)$. 
The sample variance of the $S_i$ is given by
$$
V_p := V_p(t_1,\ldots,t_p) = \frac{1}{p(p-1)} \sum_{i<j} (S_i - S_j)^2.
$$
We can exploit the previous results to compute $\BE V_p$. Since 
$$
\BE(S_i - S_j)^2 = {\rm Var} S_i + {\rm Var} S_j - 2 {\rm Cov}(S_i,S_j)  + (\BE S_i - \BE S_j)^2,
$$
we see from Corollaries~\ref{thm2} and \ref{thm3} that
\begin{eqnarray}\label{biresult}
\BE V_p(t_1,\ldots,t_p) &=& \frac{1}{p(p-1)} \, \sum_{i<j} \left\{
\theta \log\left( \frac{(e^{t_j} -e^{t_{j-1}} + e^{t_i} - e^{t_{i-1}} + 1)^2}
{(e^{t_j} -e^{t_{j-1}} + 1)(e^{t_i} - e^{t_{i-1}} + 1)}\right) \right. \nonumber\\
&& \quad + \left. \theta^2 \log^2\left( \frac{e^{t_i} - e^{t_{i-1}} + 1}{e^{t_j} -e^{t_{j-1}} + 1}
\right) \right\}.
\end{eqnarray}
If the time intervals are of equal length, say $t_i = i \tau/p$,  $i=0,1,\ldots,p$ then  (\ref{biresult}) reduces to
\begin{eqnarray}
\BE V_p(\tau) & = & \frac{\theta}{p(p-1)} \, \sum_{i<j} 
 \log\left( \frac{(\gamma^j -\gamma^{j-1} + \gamma^i  - \gamma^{i-1} + 1)^2}
{(\gamma^j - \gamma^{j-1} + 1)(\gamma^i - \gamma^{i-1} + 1)}\right) \nonumber\\
&& + \frac{\theta^2}{p(p-1)} \sum_{i<j} \log^2\left(\frac{\gamma^i - \gamma^{i-1} + 1}{\gamma^j - \gamma^{j-1} + 1}\right), \label{biresult2}
\end{eqnarray}
where $\gamma = e^{\tau/p}.$


\subsection{Logarithmically equal interval lengths}
There is one special case that results in the $S_i$ being identically distributed, namely the setting in which  $e^{t_i} = i \gamma +1$ for $i=0,1,\cdots, p$ and for some $\gamma > 0$; 
from Corollary~\ref{thm2}, the $S_i$ then have Poisson distributions with 
$$
\BE S_i = \theta \log(\gamma+1), i=1,2,\ldots,p,
$$
and, from Corollary~\ref{thm3}, covariances given by
$$
\textrm{Cov}(S_i,S_j) = \theta \log\left(\frac{(\gamma+1)^2}{2 \gamma+1}\right), i \ne j.
$$
As a consequence, the correlation between $S_i$ and $S_j$ is given by
$$
\rho := \textrm{corr}(S_i,S_j) = 2 - \log(2\gamma+1)/\log(\gamma+1), i \ne j,
$$
and, from (\ref{biresult}),
$$
\BE V_p =  \theta \log\left( \frac{2 \gamma + 1}{\gamma+1}\right).
$$
%

Given the counts $S_i, i = 1,\ldots,p$, we let $\bar{S}$ be their mean, so that
$\BE \bar{S} = \theta \log(1+\gamma)$. This suggests a Watterson-type estimator of $\theta$ given by
\begin{equation}\label{thetaS}
\theta_S = \bar{S}/\log(1+\gamma)
\end{equation}
\citep{gaw75}. The estimator is unbiased, but as $p \to \infty$,
$$
{\rm Var}\, \theta_S = {\rm Var}\, \bar{S}/\log^2(1+\gamma) = \frac{\theta}{p \log(\gamma+1)} (1 + (p-1)\rho) \to \frac{\theta \rho}{\log(\gamma+1)},
$$
so that $\theta_S$ is not a consistent estimator of $\theta$.

\subsection{Fisher's problem revisited }

Here we revisit Fisher's sampling problem \cite{raf43}. 
Several approaches have appeared in the literature, and we point the reader to \cite{SJMT22} for an overview. Here we provide an alternative setting, via the Yule process with immigration, which leads to rather transparent connections between the different views.

The setting described in \cite{SJMT22} occurs in discrete time, where successive samples of individuals of sizes $n_1, n_2, \ldots, n_p$ are taken, and the observations are $S^*_i, i=1, 2, \ldots,p$,
$S^*_i$ denoting the number of distinct types observed in the $i$th sample. The generative model is that of the CRP, a sequential model in which the distribution of the counts of family sizes follows the Ewens Sampling Formula \cite{wje72}. It is of interest to compare the two settings.

To do this, we recall first that
\begin{equation}\label{ztlaw}
\BP(Z(t) = n) = \binom{\theta+n-1}{n} e^{-\theta t} (1 - e^{-t})^n, n=0,1,2 \ldots
\end{equation}
and
\begin{equation}\label{eit}
    \BE Z(t) = \theta(e^t - 1).
\end{equation}
We could choose the time points $0 = t_0, t_1,\ldots,t_p$ in such a way that the cumulative number of individuals observed, $l_i = n_1+\cdots+n_i, i = 1,2,\ldots,p$ and  $l_0 = 0$, matches the expectation under the Yule model. To this end, we solve 
$$
 \BE Z(t_i) = \theta\left(e^{t_i}-1\right)=l_i, \ i=1,\dots,p.
 $$
to get
$$
 t_i - t_{i-1} = \log(\theta+l_i) - \log(\theta+l_{i-1}),
$$
 and 
 \begin{equation}\label{mlesoln}
 t_i=\log(\theta+l_i) -  \log \theta,
 \end{equation}
 for $i=1,\dots,p$.

We remark  that if  $\theta$ were known, and we wish to estimate the time points $t_i$ given the counts $n_i$, then (\ref{mlesoln}) gives the moment estimators of the $t_i$, and this in turn provides the maximum likelihood estimators of the $t_i$.

Henceforth we  assume that $t_i=\log(\theta+l_i) - \log \theta$, and define $\widetilde S_i = S(t_{i-1},t_i)$, $i=1,\dots,p$, and  
$$
\widetilde V_p= \widetilde V_p(n_1,\dots,n_p):=\frac{1}{p(p-1)}\sum_{i<j}(\widetilde S_i-\widetilde S_j)^2.
$$

The results of Section \ref{sect3} translate into
\begin{theorem}
For any $i\in\BN$, $\widetilde S_i$ is a Poisson random variable with mean
$$
\BE \widetilde S_i=\theta\log\left(\frac{\theta+n_i}{\theta}\right).
$$

Furthermore, for any $p,i,j\in\BN$, $i\neq j$,
$$
{\rm Cov}(\widetilde S_i,\widetilde S_j)=\theta\log\left(\frac{(\theta+n_i)(\theta+n_j)}{\theta(\theta+n_i+n_j)}\right),
$$

$$
\BE \widetilde V_p=\frac{1}{p(p-1)}\sum_{1\leq i<j\leq p}\left\lbrace\theta\log\left(\frac{(\theta+n_i+n_j)^2}{(\theta+n_i)(\theta+n_j)}\right)+\theta^2\log^2\left(\frac{\theta+n_i}{\theta+n_j}\right)\right\rbrace.
$$
\end{theorem}

\medskip
\noindent{\bf Proof. }
The results follow from substituting $t_i = \log(\theta+l_i) - \log \theta$ in Theorem~\ref{thm2}, Theorem~\ref{thm3} and (\ref{biresult}), respectively, and simplifying. \hfill\qed

\begin{corollary}
Let $n_1=n_2=\dots=n_p=n$. Then
$$
\BE \widetilde V_p =\theta\log\left(\frac{2n+\theta}{n+\theta}\right).
$$
\end{corollary}
The term on the right appears in \cite[p. 451, after (5)]{raf43} for the case $p=2$, so our model provides a unifying approach to Fisher's question.

\subsection{Asymptotic behavior}
To see the asymptotic behavior,  let $n_i=q_in$, $i=1,\dots, p$, where $q_i \in [0,1]$ satisfy $q_1+\cdots+q_p = 1$. As $n\to \infty$, we see that
$$
{\rm Cov}(\widetilde S_i,\widetilde S_j)\sim \theta\log\left(\frac{n_in_j}{n_i+n_j}\right)=\theta\log n+\theta\log\frac{q_iq_j}{q_i+q_j},
$$
and
\begin{multline*}
\BE \widetilde V_p=\frac{1}{p(p-1)}\sum_{1\leq i<j \leq p} \left\lbrace\theta\log\left(\frac{(\theta+n(q_i+q_j))^2}{(\theta+n q_i)(\theta+n q_j)}\right)+\theta^2\log^2\left(\frac{\theta+q_i n}{\theta+q_j n}\right)\right\rbrace\\[6pt]
\longrightarrow \frac{1}{p(p-1)}\sum_{1\leq i<j \leq p} \left\lbrace\theta\log\left(\frac{(q_i+q_j)^2}{q_iq_j}\right)+\theta^2\log^2\left(\frac{q_i}{q_j}\right)\right\rbrace.
\end{multline*}
The righthand formula was derived originally in \cite{bt10} by a different Poisson argument, and as a limit in the discrete case in \cite{SJMT22}.

For $n_1=\dots=n_p$, as $n\to \infty$,
$$
{\rm Cov}(\widetilde S_i,\widetilde S_j)\sim \theta\log n-\theta\log 2,
$$
and 
$$
\BE \widetilde V_p\longrightarrow \theta\log 2,
$$
as found by \cite{raf43}.

\section{How long are gaps in arrivals?}
In the covid setting it is of some interest to describe gaps in the appearance of particular variants. One setting for this is the following: Conditional on at least one family arriving in $(0,t)$, choose one of those families at random. What is the distribution of the length of time for which that family is unobservable after time $t$. Hence, we seek the distribution of the waiting time, $W_t$,  for that family to have its first birth after time $t$? 

We begin by showing that $N_t$, the number of members at time $t$ of a family randomly chosen in $(0,t)$, conditional on having at least one family, is log-series distributed, with parameter $q_t = 1 - e^{-t}$.  Since the arrival time in $(0,t)$ of a typical family is uniform on $(0,t)$, we see that
\begin{eqnarray}
\BP(N_t = j) & = & \int_0^t \frac{1}{t} \, p_j(t-u) du \nonumber \\
& = & \frac{1}{t} \int_0^t e^{-(t-u)} (1 - e^{-(t-u)})^{j-1} du \nonumber \\
& = & \frac{1}{t} \frac{q_t^j}{j} 
 =  \frac{1}{- \log(1 - q_t)} \frac{q_t^j}{j}, \label{Ndist}
\end{eqnarray}
as required.

Since each of the $N_t$ individuals in the family at time $t$ behave independently, it follows that, given $N_t = n$, $W_t$ is the minimum of $n$ independent unit exponential lifetimes, which is exponential with parameter $n$. Hence the density $f_t(s)$ of $W_t$ is
\begin{eqnarray}
f_t(s)  & = &  
\sum_{n \geq 1} \frac{1}{- \log( 1 - q_t)} \frac{q_t^ n}{n} \, \cdot ne^{-sn}, \quad s>0. \label{Wtdens}
\end{eqnarray}
This allows the moments of $W_t$ to be written down immediately:
$$
\BE W_t = \frac{1}{t} \sum_{n \geq 1} \frac{q_t^n}{n^2} = \frac{1}{t} {\rm Li}_2(q_t),\quad {\rm Var} W_t = \frac{1}{t} \sum_{n \geq 1} \frac{q_t^n}{n^3} = \frac{1}{t} {\rm Li}_3(q_t),
$$
where ${\rm Li}_n(x)$ denotes the polylogarithm function.

\section{Conditioning in the Yule process with immigration}

\cite{SJMT22} provide a discrete approach to understand Fisher's multi-sampling problem. In this paper we tackled Fisher's problem with a continuous approach.  This section discusses the connections between two approaches, focusing primarily on various versions of embedding.

To set the scene, recall that the Ewens Sampling Formula (ESF) may be realised by conditioning independent Poisson random variables on a finite \citep{gaw74a} or an infinite \citep{sl66} weighted sum. Hence we discuss a unifying approach that connects both types of conditioning relations through embedding of the CRP into the Yule process with immigration. In particular, for $\theta>0$, let
\begin{equation}\label{ESF}
\mathcal{E}_n(c_1,c_2,\dots)=\bbbone\left\lbrace\sum_{i=1}^\infty ic_i=n\right\rbrace \frac{n!}{\theta_{(n)}}\prod_{j=1}^n\left(\frac{\theta}{j}\right)^{c_j}\frac{1}{c_j!}, \ \ c_1,c_2,\dots \in \mathbb{Z}_+,
\end{equation}
for $n\in\mathbb{N}$,  where $\theta_{(n)}=\theta(\theta+1)\cdots (\theta+n-1)$ and $\theta_{(0)}=1$.  $\mathcal{E}_n$ is the ESF originally derived as the distribution of the allelic partition of a sample of $n$ genes sampled from a large population. To establish the conditioning relations, for $x\in(0,1]$, let $\pi_1(x),\pi_2(x),\dots$ be independent Poisson random variables with $\mathbb{E}\pi_i(x)=\theta x^i/i$, and let $T_n(x)=\sum_{i=1}^n i\pi_i(x)$ and $T_\infty(x)=\sum_{i=1}^\infty i\pi_i(x)$. Extending the result of \cite{sl66}, as $T_\infty(x)$ is almost surely finite for $x\in(0,1)$, we have
\begin{equation}
\mathcal{L}(\pi_1(x),\dots,\pi_n(x)\mid T_\infty(x)=n)=\mathcal{E}_n, \ \ x\in(0,1).
\end{equation}
By conditioning on $T_n(x)$ instead of the infinite sum, Watterson introduced a conditioning relation that holds for $x\in(0,1]$,  more precisely
\begin{equation}\label{relation}
\mathcal{L}(\pi_1(x),\dots,\pi_n(x)\mid T_n(x)=n)=\mathcal{E}_n, \ \ x\in(0,1].
\end{equation}

Notice from (2) that if we define $\pi_i(x)=C_i(t)$ for $t\in\mathbb{R}_+$ and $x = 1-e^{-t}\in(0,1)$, and let $Z_n(t)=\sum_{i=1}^n i C_i(t)$, by definition we have $T_n(x) = Z_n(t)$ and $T_\infty(x) = Z(t)$. Hence,  for $x\in(0,1)$ and $t\in\mathbb{R}_+$, the conditional distribution of counts of family sizes, given $Z(t)=n$ or $Z_n(t)=n$, is the ESF once more.  Note that as the $C_i(t)$ are independent Poisson random variables with $\BE C_i(t)= \theta(1-e^{-t})^i/i$,
$$
(C_1(t),C_2(t),\dots)\Rightarrow(C_1(\infty),C_2(\infty),\dots)
$$
where $C_i(\infty)=\pi_i(1)$ has a Poisson distribution with mean $ \theta/i$, for $i\in\mathbb{N}$; hence one can derive (\ref{relation}) for $x=1$,  by letting $t\to\infty$. 

\ignore{
In this section we discuss some connections between the discrete-time sampling model and its continuous-time counterpart studied here, focusing primarily on various versions of embedding. To set the scene, recall that Watterson \cite{gaw74a} showed that the Ewens Sampling Formula (ESF) may be realised by conditioning independent Poisson random variables on a weighted sum. Specifically, let $Z_1,Z_2,\ldots,Z_n$ be independent Poisson random variables with $\BE Z_i = \theta x^i/i$ for $x \in (0,1)$, and let $c_1,c_2,\ldots,c_n \in \BN$ satisfy $c_1 + 2 c_2 + \cdots + n c_n = n$. Then 
\begin{equation}\label{esfcond}
\BP(Z_1 = c_1, Z_2 = c_2,\ldots, Z_n = c_n \vert T_n = n) = \frac{n!}{\theta_{(n)}} \prod_{j=1}^n \left(\frac{\theta}{j}\right)^{c_j}\,\frac{1}{c_j!},
\end{equation}
where $x_{(n)} = x(x+1)\cdots (x+n-1), x_{(0)} = 1$ and $T_n = Z_1+2Z_2 + \cdots +n Z_n$. The distribution on the right of (\ref{esfcond}) is the ESF, originally derived as the distribution of the allelic partition of a sample of $n$ genes sampled from a large population. Notice from (\ref{Ctmeans}) that if we define $Z_i = C_i(t)$ for any fixed $t$, then we may take $x = 1 - e^{-t}$, and it follows that $T_n = Z(t)$ as defined in (\ref{ztlaw}), and hence the conditional distribution of the counts of family sizes, given $Z(t) = n$, is the ESF once more.
}

The conditioning relation provides one way to simulate observations from the ESF: simulate the $C_i(t)$, and accept the realisation if $Z(t) = n$. To make the simulation as efficient as possible, we should choose $t$ to make the probability of the conditioning event as large as possible. To do this, we note that the value of $t$ may be chosen as a function of $n$, since this choice plays no role in the conditional distribution. Maximizing $\BP(Z(t) = n)$ given in (\ref{ztlaw}) gives $t = \log((n+\theta)/\theta)$ (and so $x = x_n = n/(n+\theta)$). We note that this is the same choice of $t$ as provided in (\ref{mlesoln}). We also note that there are far more efficient ways to simulate the ESF; see \cite{abet18} for example.

The conditioning relation (\ref{relation}) is a result of the embedding of the CRP in the Yule process with immigration. Here we discuss embeddings at multiple time points. To this end, let 
$$
\tilde{\C}(n)=(\tilde{C}_1(n),\tilde{C}_2(n),\dots), \quad n\in\BN
$$
denote the counts of family sizes generated by the first $n$ arrivals in the CRP with parameter $\theta$, and let
$$
\C(t)=(C_1(t),C_2(t),\dots), \quad t\geq 0
$$
be the family size counts at time $t$ of the Yule process with immigration with birth rate $1$ and immigration rate $\theta$.  To discretise the process $\C=(\C(t), t \geq 0)$, let
$$
\Psi(n)=(\Psi_1(n),\dots,\Psi_n(n),0,0,\dots)
$$
be the family size counts of the first $n$ individuals born in the Yule process  with immigration, i.e.  $\Psi_j(n)$ gives the number of families of size $j$, considering only the first $n$ individuals. As a discrete process,  $\Psi:=(\Psi(n), n\geq 1)$   records the outcomes of jumps in the Yule process with immigration. This is a slightly different version of the jump Markov chain $J=(J(n), n\geq 1)$ used in \cite[Section 3]{st87}, in which, for each $n$, the families in $J(n)$ are also sorted in order of their appearance in the population.  In other words,  for any $n\in \BN$, one can easily obtain $\Psi(n)$ from $J(n)$ by forgetting the order (age) of the families in $J(n)$ and grouping all families of the same size together. From Theorem 2 in \cite{st87}, $J$ is independent of $Z$, then as a functional of $J$, so is $\Psi$. The independence comes as a result of the fact that at each jump time, one new individual will be added to the existing population, no matter if the new individual arrives as a result of an immigration (new family) or a birth. The connection between $\tilde{\C}=(\tilde{\C}(n), n\in \BN)$ and $\C=(\C(t)), t \geq 0)$  is given in the next theorem.

\begin{theorem}
For any $p\in \BN$, $0<t_1<\dots<t_p$, and $0\leq l_1\leq l_2\leq \dots\leq l_p$, $l_i\in\BZ_+$, and $u_1,u_2,\dots,u_p\in \BZ_+^\BN$, we have
\begin{multline}\label{thmcond}
\BP(\C(t_1)=u_1,\dots,\C(t_p)=u_p\mid Z(t_1)=l_1,\dots,Z(t_p)=l_p)\\
=\BP(\tilde{\C}(l_1)=u_1,\dots,\tilde{\C}(l_p)=u_p).
\end{multline}
\end{theorem}

\medskip
\noindent{\bf Proof. }
First note that if $Z(t_1)=l_1,\dots,Z(t_p)=l_p$, then $\Psi(l_i)=\C(t_i)$, for $i=1,\dots, p$. On the other hand, (3.1) in \cite{st87} and the connection between $\Psi$ and $J$ shows that $(\Psi(i), i\geq 1)\sim (\tilde{\C}(i), i\geq 1)$ in distribution. Also, from Theorem 2 in \cite{st87} $J$ and $Z$, and hence $\Psi$ and $Z$ are independent processes. Thus we can write
\begin{eqnarray*}
\frac{\BP(\C(t_i)=u_i,Z(t_i)=l_i;i=1,\dots,p)}{\BP(Z(t_i)=l_i;i=1,\dots,p)}&=&\frac{\BP(\Psi(l_i)=u_i,Z(t_i)=l_i;i=1,\dots,p)}{\BP(Z(t_i)=l_i;i=1,\dots,p)}\\[4pt]
&=&\BP(\Psi(l_1)=u_1,\dots,\Psi(l_p)=u_p)\\[4pt]
&=&\BP(\tilde \C(l_1)=u_1,\dots,\tilde \C(l_p)=u_p).
\end{eqnarray*}
\hfill\qed

To better connect the sequential multi-sampling theory of the Yule process with immigration to its discrete-time counterpart, consider a population of size $n_1+n_2+\cdots+n_p$, sampled from a single run of a CRP.  \cite{SJMT22} study the pairwise correlation and sample variance of $S_1^*,S_2^*, \cdots,S_p^*$, the number of types (or species) appearing in the first $n_1$ arrivals, the second $n_2$ arrivals, \ldots, and the last $n_p$ arrivals of the CRP sample.  It is now straightforward from (\ref{thmcond}) that $(S_1^*,S_2^*, \cdots,S_p^*)$ is in distribution the same as $(S(t_0,t_1),\cdots,S(t_{p-1},t_p))$, conditional on observing,  in the latter, exactly $n_i$ individuals in $(t_{i-1},t_i)$, for $i=1,\cdots,p$. As mentioned in the discussion after (\ref{mlesoln}), letting $t_i=\log((\theta+l_i)/\theta)$ for $i=1,\cdots,p$, maximizes the chance of observing $n_i$ new individuals (i.e. $n_i$ births) in $(t_{i-1},t_i)$, and under this assumption, the asymptotic behavior of the sample variance of $\tilde{S}_1,\tilde{S}_2,\cdots,\tilde{S}_p$ and that of the sample variance of $S_1^*,S_2^*, \cdots,S_p^*$ coincide. In this case, in addition to providing a relatively simpler way to calculate things, the Yule process with immigration allows the sequential samples $\tilde{S}_1,\cdots,\tilde{S}_p$ to have a random number of individuals, and hence is more appropriate for population models with random sample sizes in which the individuals arrive one by one at random times.

\section*{Acknowledgements}
PHdS and ST were supported in part by National Science Foundation grant DMS2030562. 

\bibliography{Arxiv_Fisher_BI_11-28-2022}
\bibliographystyle{plainnat}


\end{document}